\documentclass[a4j,11pt]{article}
\usepackage{amsthm}
\usepackage{amsmath}
\usepackage{amssymb}
\usepackage{amsfonts}
\usepackage{wrapfig}
\usepackage[pdftex]{graphicx,color}

\theoremstyle{plain}

\newtheorem{thm}{Theorem}
\newtheorem{remark}{Remark}
\newtheorem{col}{Corollary }
\newcommand{\argmax}{\mathop{\rm arg~max}\limits}

\title{Non-area-specific adjustment factor for second-order efficient empirical Bayes confidence interval}
\author{Masayo Yoshimori Hirose \\{The Institute of Statistical Mathematics} }
\date{\hfill}

\begin{document}

\maketitle
\begin{abstract}
An empirical Bayes confidence interval has high user demand in many applications. In particular, the second-order empirical Bayes confidence interval, the coverage error of which is of the third order for a large number of areas, $m$, is widely used in small area estimation when the sample size within each area is not large enough to make reliable direct estimates according to a design-based approach. Yoshimori and Lahiri (2014a) proposed a new type of confidence interval, called the {\it second-order efficient empirical Bayes confidence interval}, with a length less than that of the direct confidence estimated according to the design-based approach. However, this interval still has some disadvantages: (i) it is hard to use when at least one leverage value is high; (ii) many iterations tend to be required to obtain the estimators of one global model variance parameter as the number of areas, $m$, increases, due to the area-specific adjustment factor.
To prevent such issues, this study proposes a more efficient confidence interval to allow for high leverage and reduce the number of iterations for large $m$, by adopting a non-area-specific adjustment factor and coordinate the measure of uncertainty of the empirical Bayes estimator, maintaining the existing desired properties. Moreover, we present two simulation studies to show the efficiency of this confidence interval.

\end{abstract}


\section{Introduction}
There has been increasing demand for reliable statistics of government fund allocations, social services planning, etc., in smaller geographic areas and sub-populations, where large samples are not available. Because of the limited number of observations within each area or domain, a direct estimator constructed according to the design-based approach only from information within each area or domain, is not reliable.
The empirical Bayes estimator and empirical best linear unbiased predictor (EBLUP) help make efficient inferences by borrowing information from other areas via model-based approaches to small area estimation.
Fay and Herriot (1979) first applied this model-based approach to Census data through a specific Bayesian model. The model, called the Fay--Herriot model, has been widely used in practice.
For $i=1,\ldots,m,$
\begin{align}
\noindent {\rm Level \ 1:}& y_{i}\mid\theta_i \stackrel{ind.}{\sim} N(\theta_i,D_i);\notag\\
\noindent {\rm Level\ 2:}& \theta_i \stackrel{ind.}{\sim} N(x^{\prime}_i\beta,A).\label{FH}
\end{align}
In the above model, level 1 is used to take into account the sampling distribution of the direct estimator $y_i$ for small area $i$.
A true mean for small area $i$, $\theta_i$, is linked to provide the auxiliary variables $x_i=(x_{i1},\cdots,x_{ip})$ in a level-2 linking model. In practice, the coefficient $p$-vector $\beta$ and the model variance parameter $A$ in the linking model are unknown, and we need to estimate them from the observed data. The assumption of a known $D_i$ often follows from the asymptotic variances of the transformed direct estimates (Efron and Morris 1975) or from empirical variance modeling (Fay and Herriot 1979).
This model can be viewed as the following linear mixed model:
$$y_i=\theta_i+e_i=x_i^{\prime}\beta+u_i+e_i, \ i=1,\ldots,m,$$
where $u_i$ and $e_i$ are independent of the normality assumption $u_i \stackrel{i.i.d.}{\sim}N(0,A)$ and $e_i\stackrel{ind.}{\sim}N(0,D_i)$.

Let $M_i$ define the mean squared error (MSE) $E[(\hat{\theta}_i-\theta_i)^2]$ of the predictor $\hat{\theta}_i$ of a small area mean $\theta_i$, where the expectation is on the joint distribution of $y$ and $\theta$ under the Fay--Herriot model (\ref{FH}) with $y=(y_1,\ldots,y_m)^{\prime}$ and $\theta=(\theta_1,\ldots,\theta_m)^{\prime}$.

\indent The Bayes estimator of $\theta_i$ is consistent with the best predictor (BP) in this model, with the minimum MSE among all $\hat{\theta}_i$. It is given by
$$\hat{\theta}_i^{BP}=(1-B_i)y_i+B_i x^{\prime}_i\beta,$$
where $B_i =\frac{D_i}{A+D_i}$ is called the shrinkage factor toward $x^{\prime}_i\beta$ from the direct estimate $y_i$.

If $\beta$ is unknown, the best linear unbiased predictor (BLUP), in which $\beta$ of $\hat{\theta}_i^{BP}$ is replaced by $\tilde{\beta}$, minimizes the MSE among all linear unbiased predictors of $\theta_i$, as follows:
$$\hat{\theta}_i^{BLUP}=(1-B_i)y_i+B_i x^{\prime}_i\tilde{\beta},$$
where the weighted least-square estimator of $\beta$, $\tilde{\beta}=\tilde{\beta}(A)=( X^{\prime}V^{-1}X)^{-1} X^{\prime}V^{-1} y$, $y=(y_1\ldots,y_m)^{\prime}$, $X=(x_1,\ldots,x_m)^{\prime}$ and $V=diag(A+D_1,\ldots,A+D_m)$. \\

\indent From the fact that both $\beta$ and $A$ are practically unknown, the empirical best linear unbiased predictor (EBLUP) $\hat{\theta}_i^{EB}$ is widely used for small area inference, where the unknown model variance parameter $A$ in $\hat{\theta}^{BLUP}_i$ is replaced by its consistent estimator $\hat{A}$:
$$\hat{\theta}_i^{EB}=(1-\hat{B}_i)y_i+\hat{B}_i x^{\prime}_i\hat{\beta},$$
where $\hat{B}_i=\frac{D_i}{\hat{A}+D_i}$ and $\hat{\beta}(\hat{A})=\hat{\beta}=\tilde{\beta}(\hat{A})$ and the consistent estimator $\hat{A}$ for large $m$ is even translation invariant for all $y_i$ and $\beta$.
To estimate the model variance parameter $A$, the methods of moments estimator (see Fay and Herriot 1979; Prasad and Rao 1990) and standard maximum likelihood estimators, such as profile maximum likelihood (ML) estimator and residual maximum likelihood (REML) estimator are utilized. In particular, the REML estimator of $A$ is widely used in terms of higher-order asymptotic properties for large $m$ under some mild regularity conditions.
Hereafter, we indicate the REML estimator as $\hat{A}_{RE}$, obtained from $$\hat{A}_{RE}=\argmax_{0\leq A<\infty} L_{RE}(A|y),$$
where the residual likelihood function $L_{RE}(A|y)=|X^{\prime}V^{-1}X|^{-1/2}|V|^{-1/2}\exp\{ -y^{\prime}Py/2\}$ and $P=V^{-1}-V^{-1}X( X^{\prime}V^{-1}X)^{-1}X^{\prime}V^{-1}$.

This study focuses on the confidence interval for $\theta_i$, used widely in small area estimation as well as point estimation.
Let $I_{i}$ denote the general form of the confidence interval as follows:
 \begin{align}
I_{i}: {\xi}_i \pm q_{i}s_i, \label{CI.form}
\end{align}
where $\xi_i$ and $s_i$ are, respectively, a predictor of $\theta_i$ and a measure of uncertainty of $\xi_i$. $q_i$ is adopted as an adequate percentile point to get closer to the nominal coverage level $1-\alpha$.

We call the $100(1-\alpha)$\% confidence interval for $\theta_i$ if the coverage probability is consistent with the nominal coverage, that is, $P[\theta_i\in I_i|\beta,A]=1-\alpha$ holds exactly for a fixed $\beta$ and $A$ with a probability measure $P$, according to the Fay--Herriot model.
We introduce several intervals traditionally used for small area estimation.
The simplest confidence interval, called the direct confidence interval and denoted by $I_i^{D}$, is constructed with the direct estimate $y_i$, the $z$ value $z_{\alpha/2}$ of the upper $100(1-\alpha)$\%, and the sampling variance $D_i$, substituted for  $\xi_i$, $q_i$, and $s_i^2$, respectively.
This yields a coverage probability of exactly $1-\alpha$. However, this interval could yield an excessively long length to make any reasonable conclusion when $D_i$ is large.
In contrast, Cox (1975) suggested the empirical Bayes confidence interval for $\theta_i$ in balanced case, choosing ${\xi}_i=\hat{\theta}_i^{EB}(\hat{A})$, $q_i=z_{\alpha/2}$, $s_i^2={g_{1i}(\hat{A})}$ and plugging the ANOVA estimator $\hat{A}_{ANOVA}$ into $\hat{A}$, denoted by $I_i^{Cox}(\hat{A}_{ANOVA})$, where $g_{1i}(A)=\frac{AD_i}{A+D_i}$.
Although successful, because its length is not greater than that of the direct confidence interval $I^{D}_i$,
it is not accurate enough in most small area applications because the coverage error is of the order of $O(m^{-1})$ for large $m$.

Similar to the Cox interval, the traditional empirical Bayes confidence interval $I_{i}^{T}$ is also a common method suggested in Prasad and Rao (1990), in which ${\xi}_i$, $q_i$, and $s_i$ in (\ref{CI.form}) are replaced by $\hat{\theta}_i^{EB}$, $ z_{\alpha/2}$, and $s_i^2={\hat{M}_i(\hat{A}_{PR})}$, along with the PR estimator using Henderson\rq{}s method 3 (Prasad and Rao 1990) of $A$, where the second-order unbiased estimator $\hat{M}_i$ of $M_i$ is such that $E[\hat{M}_i-M_i]=o(m^{-1})$.

Recently, users have frequently employed the residual maximum likelihood estimator $\hat{A}_{RE}$ instead of the ANOVA-type estimators for two such intervals. We denote these as $I_{i}^{Cox}\equiv I_{i}^{Cox}(\hat{A}_{RE})$ and $I_{i}^{T}\equiv I_{i}^{T}(\hat{A}_{RE})$, with the REML estimator $\hat{A}_{RE}$ used for each empirical Bayes confidence interval.
The length of $I_i^{T}$ obtained is less than that of the direct confidence interval for large $m$, but the coverage error is of the order of $O(m^{-1})$ as with the Cox interval $I_{i}^{Cox}$.

To offset this disadvantage, Datta et al. (2002) and Sasase and Kubokawa (2005) suggested calibrating $z_{\alpha/2}$ to $q_i$ for the Cox confidence interval $I_{i}^{Cox}$ to reduce the coverage error to the order of $O(m^{-3/2})$, for a balanced case ($D_i=D$ for all $i$) in the Fay--Herriot model and for an unbalanced case in the nested error regression model.
Diao et al. (2014) revealed the relationship between the percentile point $q_i$ and the measure of uncertainty of $\hat{\theta}_i^{EB}$, $s_i$, in constructing the second-order confidence interval based on the existing variance estimation method. From this result, if REML is adopted as the estimator of $A$, it assigns ${\xi}_i=\hat{\theta}_i^{EB}(\hat{A}_{RE})$, $q_i=z^*$ with $z^*$ as stated in Corollary 1 in Diao et al. (2014) in constructing second-order confidence interval $I_{i}^{CT}$, $s_i^2={\hat{M}_i(\hat{A}_{RE})=g_{1i}(\hat{A}_{RE})+g_{2i}(\hat{A}_{RE})+2g_{3i}(\hat{A}_{RE})}$ in (\ref{CI.form}), where $g_{2i}(A)=B_i^2x_i^{\prime}(X^{\prime}V^{-1}X)^{-1}x_i$ and $g_{3i}(A)=2B_i^2/\{(A+D_i)tr(V^{-2})\}$. 
The second-order confidence interval can also be constructed through a simulation-based method.
Hall and Maiti (2006), Chatterjee et al. (2008), and Li and Lahiri (2010) achieved an adequate $q_i$ by a bootstrap method. Yoshimori (2015) found, by a simulation study based on the Fay--Herriot model, that the interval length of Chatterjee et al.(2008) tends to be less than that of Hall and Maiti (2006).
Yoshimori and Lahiri (2014 a) suggested a second-order efficient empirical Bayes confidence interval via a new area-specific adjustment factor with fixed $s_i$ as $\sqrt{g_{1i}(A)}$, which first achieves three desired properties simultaneously: (i) the coverage error is of the order $O(m^{-3/2})$, (ii) the length is always less than that of the direct confidence interval, and (iii) it does not rely on a simulation-based method such as the bootstrap method. Hereafter, we call this the YL confidence interval, which is denoted by $I_{i}^{YL}$ and written as follows:
\begin{align}
\hat \theta^{EB}(\hat{A}_{i,YL})\pm z_{\alpha/2}\sqrt{g_{1i}(\hat A_{i,YL}) } \label{YL}
\end{align}
where $\hat{A}_{i, YL}$ yields the maximum value of their adjusted residual likelihood, $L_{ad}^{YL}$, which is equal to the residual likelihood $L_{RE}(A|y)$ multiplied by a specific adjustment factor, $\tilde{L}_{i,ad}^{YL}(A)$; that is, $\hat{A}_{i,YL}=\argmax_{A>0} \tilde{L}_{i,ad}^{YL}(A)L_{RE}(A|y)$, where $\tilde{L}_{i,ad}^{YL}(A)$ is as shown in Yoshimori and Lahiri (2014a).
For the concept of the adjusted likelihood method, refer to Lahiri and Li (2009). 

The authors emphasized that the YL interval method would be an alternative to the parametric bootstrap interval proposed in Li and Lahiri (2010), in spite of the fact that YL does not rely on a simulation-based method. From the result, it would be reliable in small area estimation, especially for developing countries that are not in favor of using simulation-based methods.
Nevertheless, their required condition for the existence of $\hat{A}_{i,YL}$ could be getting stronger not for large $m$ when at least one leverage value is high due to the following condition in practice.
\begin{align}
m>\frac{4+p}{1-h_i},\label{cond.YL}
\end{align}
where $h_i$ is the leverage $x_i^{\prime}(X^{\prime}X)^{-1}x_i$.
This is the condition described in Remark 3 of Yoshimori and Lahiri (2014a) to ensure that $\hat{A}_{i,YL}$ exists.
For example, this condition does not hold for the baseball batting average data in using previous seasonal batting average as one covariate, shown in Gelman et al. (1995).

Moreover, the YL interval method might cause confusion in two ways: (a) some people might not understand why $m$ estimates are needed for one global parameter $A$; (b) with no parallel computations for large $m$, considerable time might be needed for the number of iterations required to obtain area-specific $m$ estimates of $A$, with the likelihood method used only for the global parameter $A$.

To address these problems, this study proposes a pioneering, more reliable confidence interval, see in (\ref{NAS}), satisfying the following five desired properties, by providing a new non-area-specific (NAS) adjustment factor and coordinate the measure of uncertainty of EBLUP based on the NAS method with fixed $q_i=z_{\alpha/2}$:

{\bf Desired properties}
\begin{description}
\item [(i)] The coverage error is of the order of $O(m^{-3/2})$;
\item [(ii)] The length is always less than that of the direct confidence interval;
\item [(iii)] It does not rely on a simulation-based method, such as the bootstrap method;
\item [(iv)] It does not require calculations to obtain the estimator of $A$ for all $m$ areas, unlike the YL method (2014a);
\item [(v)] It has a milder condition for the existence of the estimator of $A$ than the (\ref{cond.YL}).
\end{description}

Note that Property (i) does not hold for $I^{Cox}$ and $I_i^{T}$, while Property (ii) does not hold for $I_i^{T}$ and $I_i^{CT}$. In addition, the simulation-based method does not satisfy Properties (ii)--(iii).

We also show some simulation results for comparison with several existing intervals. 

Now, we prepare the regularity conditions for the introduction of some theorems mentioned in the next section, which correspond to R2--R4 in Yoshimori and Lahiri (2014a).

{\bf Regularity conditions}
\begin{description}
\item [R1]$rank(X) = p$ is bounded for large $m$;
\item [R2]The elements of $X$ are uniformly bounded, implying $\sup_{i \geq 1} h_i = O(m^{-1})$; 
\item[R3]$0 < \inf_{i\geq1}D_i \leq \sup_{i\geq1}D_i <\infty$, $A \in (0,\infty)$;
\end{description}
Hereafter, we constrain the class of the adjustment factor $\tilde{L}_{i,ad}(A)$ with the following conditions, which correspond to conditions R1 and R5 in Yoshimori and Lahiri (2014a).
\begin{description}
\item [R4]The logarithm of the adjustment factor $\tilde{l}_{ad}(A)$ [or $\tilde{l}_{i,ad}(A)$] is free of $y$;
it is also five times continuously differentiable with respect to $A$. Moreover, it is bounded for large $m$;
\item [R5]$|\hat{A}_{i}|<C_+m^{\lambda}$, where $C_+$ is a generic positive constant and $\lambda$ is a small positive constant, where $\hat{A}_i=\argmax_{A>0} \tilde{L}_{i,ad}(A)L_{RE}(A|y)$ for an adjustment factor $\tilde{L}_{i,ad}(A)$.
\end{description}

\section{Second-order efficient confidence interval based on a non-area-specific adjustment factor}

We consider the following confidence interval class with general adjusted residual maximum likelihood estimator $\hat{A}_i$ and $c_i^*(\hat{A}_i,z)$: 
\begin{align}
\label{class1}
I_{i}^{*}:\theta_i^{EB}(\hat{A}_i) \pm z_{\alpha/2}\sqrt{g_{1i}(\hat{A}_i) +g_{2i}(\hat{A}_i)+\hat{c}_i^{*}(\hat{A}_i,z)g_{3i}(\hat{A}_i)},
\end{align}
where the function $c_i^*=c_i^*(A,z)$ is defined as $\{s_i^2-(g_{1i}(A)+g_{2i}(A))\}/g_{3i}(A)$ with $s_i$ as given in (\ref{CI.form}). 
Note that $c_i^*$ is set as $-g_{2i}(A)/g_{3i}(A)$ in the YL method and the class $c_i^*$ is restricted to 
$$c_i^*>-\frac{(A+D_i)^2tr[V^{-2}]}{2D_i}\left\{A+\frac{D_i}{A+D_i}x_{i}^{\prime}(X^{\prime}V^{-1}X)^{-1}x_{i} \right\}$$
since $s_i^2>0$.

In this study, we replace $c_i$ in Diao et al.(2014) by $c_{i}^*({A},z)g_{3i}({A}_i)$ as $g_{3i}(A)$ is related to the uncertainty measure of EBLUP, ${M}_i(\hat{\theta}_i^{EB})$.

We first construct the following theorem to show the relationship between the adjustment factor $\tilde{L}_{i,ad}(A)$ and the uncertainty measure of $\hat{\theta}^{EB}_i$, that is, $s_i^2=g_{1i}(\hat{A}_i) +g_{2i}(\hat{A}_i)+c_i^{*}(\hat{A}_i,z)g_{3i}(\hat{A}_i)$:
\begin{thm}
\label{t1}
{\rm 
Under the regularity conditions, we use the following equation to construct a second-order empirical Bayes confidence interval for large $m$:
\begin{align}
\tilde{l}_{i,ad}^{(1)}(A)&=\frac{7-z^2-4c_i^{*}(A,z)}{4(A+D_i)}+\frac{(1+z^{2})}{4A}+O(m^{-1/2}), \label{ad.diff2}
\end{align}
where the preassigned $z=z_{\alpha/2}$ and $\tilde{l}_{i,ad}^{(1)}$ is the first derivative of the logarithm of the adjustment factor $\tilde{L}_{i,ad}(A)$ with respect to $A$. The proof is given in Appendix A.
}
\end{thm}

Theorem \ref{t1} ensures several corollaries as follows:
\begin{col}
{\rm 
\begin{enumerate}
\item [(a)] When we adopt an REML estimator as $\hat{A}$, a suitable $c_i^*$ can be derived with the right-hand side of (\ref{ad.diff2}) set to zero, that is, $\tilde{l}_{i,ad}^{(1)}(A)=0$. 
Unfortunately, it does not satisfy the desired property (ii), as $c_i^*$ is derived as $\hat{c}_i^*=\left[2+\frac{(1+z^2)D_i}{4\hat{A}}\right ]$, depending on $A$, which causes a substantial increase in the variability of the length as $\hat{A}$ decreases.
\item [(b)] When $c_i^*=0$, we use $g_{1i}(A)+g_{2i}(A)$ as $s_i^2$. The following confidence interval can then be constructed:
\begin{align}
I_{i}^{(c)}: \hat{\theta}_i^{EB}(\hat{A}_i^{(c)})\pm z_{\alpha/2}\sqrt{g_{1i}(\hat{A}_i^{(c)})+ g_{2i}(\hat{A}_i^{(c)})},\label{sugg2}
\end{align}
where $\hat{A}_i^{(c)}=\arg\max_{A>0}\tilde{L}_{i,ad}^{(c)}(A)L_{RE}(A|y)$, with $$\tilde{L}_{i,ad}^{(c)}(A)=A^{\frac{1+z^2}{4}}(A+D_i)^{\frac{7-z^2}{4}}$$ and $g_{1i}(A)+g_{2i}(A)$, is the mean squared error of BLUP, $M[\hat{\theta}_{i}^{BLUP}]$.
This interval also satisfies properties (i)--(iii), since $g_{1i}(A)+g_{2i}(A)<D_i$ holds where $A>0$, the proof of which is given in Appendix B, but property (iv) is sacrificed except in a balanced case. 
\item[(c)] If $c_i^{*}$ is set to 2 for all $i$ in the interval class (\ref{class1}), as with the traditional empirical Bayes confidence interval $I_i^{T}$, we obtain the adjustment factor as $\tilde{L}_{i,ad}(A)=\frac{(1+z^2)D_i}{4A(A+D_i)}$ for the second-order confidence interval, which cannot be a non-area-specific adjustment factor.
\end{enumerate}
}
\end{col}

From Corollary 1(a), we consider the class of $c_i^*$ such that $c_i^{*}$ is independent of $A$, that is $c_i^{*}(A,z)=c_i^*(z)$, to avoid a substantial increase in the variability of the length, unlike in the YL method.

Then, to satisfy property (iv) as well, our first goal is to find an adequate non-area-specific adjustment factor $\tilde{L}_{i,ad}(A)=\tilde{L}_{ad}(A)$ such that it is not free from $A$, in terms of Corollary 1(a).
Remember that satisfying property (iv) does imply that iterative calculations are not required to estimate $A$ for the whole $i$ area in using a likelihood method.

The result shows that there is only one non-area-specific term $\frac{(1+z^{2})}{4A}$ in equation (\ref{ad.diff2}), up to the order $O(m^{-1/2})$. 
If we eliminate the remaining term, we achieve our first goal successfully; that is,
$$\frac{7-z^2-4c_i^*(z)}{4(A+D_i)} =0.$$
It follows that $c_i^*(z)$ no longer depends on $i$ such that $c_i^*=c^*=\frac{(7-z^2)}{4}$ for all $i$. 
This implies that we achieve our first goal by combining $s_i^2=g_{1i}(A)+g_{2i}(A)+\frac{(7-z^2)}{4}g_{3i}(A)$ with the non-area-specific adjustment factor $$\tilde{L}_{ad}(A)= A^{\frac{(1+z^2)}{4}},$$ obtained from the differential equation
$\tilde{l}_{ad}^{(1)}=\frac{(1+z^2)}{4A}$. 
Hereafter, we denote this non-area-specific adjustment factor as $\tilde{L}_{NAS}(A)$.

According to this result, we construct the second-order confidence interval $I_i^{NAS0}$ by combining the non-area-specific adjustment factor $\tilde{L}_{ad}$ with $c^*(z)$ for all $i$:
\begin{align}
I_i^{NAS0}: \hat \theta^{EB}(\hat{A}^{NAS})\pm z_{\alpha/2}\sqrt{g_{1i}(\hat A^{NAS}) +g_{2i}(\hat A^{NAS})+ \frac{(7-z^2)}{4}g_{3i}(\hat A^{NAS})} \label{sugg1}
\end{align}
and
$\hat{A}^{NAS}=\arg\max_{A >0} \tilde{L}_{NAS}(A)L_{RE}(A|y)$ with non-area specific adjustment factor $\tilde{L}_{NAS}(A)$.

\begin{thm}
{\rm
\label{t.2}
Under the regularity conditions R1-R3 and R5, 
\begin{enumerate}
\item[(1)] there exists at least one solution $\hat{A}^{NAS}$ for $A>0$ with the condition $m>p+\frac{1+z^2}{2}$. 
\item[(2)] In balanced case, there is a unique solution $\hat{A}^{NAS}$ for $A>0$.
\end{enumerate}
}
\end{thm}

The proof is given in Appendix A. 

Theorem \ref{t.2} ensures that this existence condition for the estimator of $A$ is no longer dependent on leverage $h_i$ and is milder than (\ref{cond.YL}) when $z_{\alpha/2}^2<7$. 
This implies that our interval $I_{i}^{NAS0}$ has property (v) as well.

In the discussion so far, $I_{i}^{NAS0}$ has shown properties (i), (iii), (iv), and (v). Does $I_{i}^{NAS0}$ then have property (ii) as well?
In general, its length is always less than $2z_{\alpha/2}\sqrt{D_i[1+\frac{7-z^2}{2}]}$, but not less than that of the direct confidence interval. Even so, the length does not tend to be inflated significantly since $c^*(z)=\frac{7-z^2}{4}$ does not depend on $\hat{A}_i$.

\begin{remark}
{\rm 
We can also consider the possibility of another non-area-specific factor not dependent on $z$, $\tilde{L}_{ad}(A)=A^{1/4}$,
with $c_i^*=\left[\frac{7}{4}+\frac{z^2D_i}{4\hat{A}}\right ]$ such that
\begin{align*}
\hat \theta^{EB}(\hat{A})\pm z_{\alpha/2}\sqrt{g_{1i}(\hat{A}) +g_{2i}(\hat A)+\left[\frac{7}{4}+\frac{z^2D_i}{4\hat{A}}\right ]g_{3i}(\hat{A})}
\end{align*}
where $\hat{A}=\arg\max_{A>0}A^{1/4}L_{RE}(A|y)$.

However, this could increase the variability of the length since the estimator of $c_i^*$ is dependent on $\hat{A}$.
}
\end{remark}

We next set the second goal: construct an interval such that all desired properties described in the previous section are satisfied.
To achieve the goal, we suggest the following second-order efficient confidence interval using $I_{i}^{(c)}$ given in Corollary 1 (b):
\begin{eqnarray}
\label{NAS}
I_i^{NAS}: \left\{ \begin{array}{ll}
I_i^{NAS0} & i \in S\equiv\{i \mid{g_{1i}(\hat{A}^{NAS})+ g_{2i}(\hat{A}^{NAS})+\frac{7-z^2}{4}g_{3i}(\hat{A}^{NAS})}<{D_i}\} \\
I_{i}^{(c)} & {\rm otherwise \ while \ the \ regularity \ conditions \ hold.}\\
\end{array} \right.
\end{eqnarray}
Thus, we achieve our second goal successfully.

Our suggested estimators $\hat{A}^{NAS}$ and $\hat{A}_{i}^{(c)}$ have the following properties as well, which were used in $I_i^{NAS}$.

\begin{thm}
{\rm
Under the regularity conditions R1-R3 and R5, we have, for large $m$,
\label{t3}
\begin{enumerate}
\item [(1)] $E[\hat{A}_i^{*}-A]=\frac{2}{tr[V^{-2}]}\tilde{l}_{i,ad*}^{(1)}+o(m^{-1}) $; 
\item [(2)] $E[(\hat{A}_i^{*}-A)^2]= \frac{2}{tr[V^{-2}]}+o(m^{-1})$;
\item [(3)] There exists at least one solution $\hat{A}_i^{(c)}$ on $A>0$ for $i$ with the condition $m>p+4$,

\end{enumerate}
where $\hat{A}_i^{*} \in \{\hat{A}^{NAS}, \ \hat{A}_{i}^{(c)}\}$ and $\tilde{l}_{i,ad*}^{(1)} \in \{\frac{1+z^2}{4A}$, \ $\frac{7-z^2}{4(A+D_i)}+\frac{1+z^2}{4A}$\} for each $\hat{A}^{NAS}$, $\hat{A}_{i}^{(c)}$.
}
\end{thm}
Proofs (1)-(2) follow directly from Theorem 1 in Yoshimori and Lahiri (2014b) or the proof of the corollary to Theorem 4 in Yoshimori and Lahiri (2014a). Part (3) implies that $I_i^{(c)}$ has property (v) as well and 
that proof is deferred to Appendix A.

\section{Simulation Study}
\subsection{Simulation design}
In this section, two finite sample simulation studies are implemented to investigate the performances of non-simulation-based confidence intervals through Monte Carlo simulation under the Fay--Herriot model (\ref{FH}).
The first simulation study was considered a balanced case such that the leverages of all areas satisfy the required condition (\ref{cond.YL}) for the confidence interval $I_i^{YL}$, whereas an unbalanced case with one area does not satisfy the condition (\ref{cond.YL}) for the second study.
We set $\beta=0$ without loss of generality through these simulation studies. For each study, we generated $10^4$ datasets from model (\ref{FH}). When the REML method obtained zero estimates, we truncated it to 0.01.

\subsection{Study 1: balanced case in which condition (\ref{cond.YL}) holds for all areas}
For the first simulation study, we considered a situation in which the number of areas was $m=15$ and the dimensions of $\beta$ were $p=2$ and $x_{i1}\beta_1=\mu$ for area $i$. Additionally, we generated $x_{i2}$ independently from the uniform distribution $U(0,1)$ once and then treated it as fixed. Then, all leverages satisfied condition (\ref{cond.YL}) since the maximum was 0.23. In order to investigate the effect of the shrinkage factor $B_i=B$ in a balanced case, we considered three $B$ values, 0.5, 0.7, and 0.9, 
defining the sampling variances such that $D_i=D=1$ and changing the value of $A$. 
The competitors are the Cox confidence interval with the REML method; the traditional empirical Bayes confidence interval with the REML method, proposed in Prasad and Rao (1990); the Cox-type confidence interval with $\hat{A}_{i,YL}$ given in (\ref{YL}); the calibrated traditional confidence interval with the REML method, proposed in Diao et al. (2014); the second-order efficient confidence interval based on a non-area-specific adjustment factor, proposed in this study; and the direct confidence interval. 
Let them be denoted by {\it Cox.Re, T.Re, Cox.YL, CT.Re, NAS, and Direct}, respectively.

\begin{table}[h]
\centering
\small
\caption{Simulated coverage probabilities and average length (in parentheses) in a balanced case such that condition (\ref{cond.YL}) holds with 95\% nominal coverage}
\label{balance}
\begin{tabular}{c|c|ccccccc}
  \hline
B & Leverage & Cox.Re &T.Re& Cox.YL & NAS & CT.Re & Direct \\
  \hline
0.50 & 0.07 & 82.81 & 94.05 & 95.87 & 96.24 & 98.08 & 94.87 \\ 
& & ( 2.44 ) & ( 3.07 ) & ( 3.27 ) & ( 3.23 ) & ( 6.25 ) & ( 3.92 ) \\ \cline{2-8}
& 0.23 & 80.75 & 95.4 & 95.57 & 96.18 & 97.9 & 95.49 \\ 
&  & ( 2.44 ) & ( 3.31 ) & ( 3.48 ) & ( 3.38 ) & ( 3.97 ) & ( 3.92 ) \\  \hline
0.70 & 0.07 & 72.19 & 96.63 & 96.74 & 97.66 & 99.68 & 95.37 \\ 
&  & ( 1.76 ) & ( 2.77 ) & ( 3.07 ) & ( 3.04 ) & ( 11.71 ) & ( 3.92 ) \\ \cline{2-8}
& 0.23 & 68.25 & 96.6 & 95.78 & 97.08 & 99.37 & 95.02 \\ 
&  & ( 1.76 ) & ( 3.1 ) & ( 3.33 ) & ( 3.24 ) & ( 4.6 ) & ( 3.92 ) \\  \hline
0.90 & 0.07 & 66.14 & 99.47 & 97.7 & 98.82 & 99.88 & 94.65 \\ 
&  & ( 1.21 ) & ( 2.56 ) & ( 2.91 ) & ( 2.89 ) & ( 17.95 ) & ( 3.92 ) \\ \cline{2-8}
& 0.23 & 58.14 & 98.84 & 96.17 & 98.2 & 99.82 & 94.94 \\ 
&  & ( 1.21 ) & ( 2.96 ) & ( 3.21 ) & ( 3.13 ) & ( 5.27 ) & ( 3.92 ) \\ 
   \hline
\end{tabular}
\end{table}

Table \ref{balance} shows the simulated coverage probability and average length within parentheses with 95\% nominal coverage for each $B$ in combination with minimum and maximum leverage values of 0.07 and 0.23, respectively.
 As shown in this table, even in the case of minimum values (Leverage 0.07 and $B$=0.5), the results for $I_i^{Cox}$and $I_i^{T}$ indicate under-coverage although the length of the direct confidence interval decreases in comparison with the intervals based on other methods. In practice, a severe problem occurs from under-coverage, which might be affected by the absence of Property (i) of $I_i^{Cox}$and $I_i^{T}$.
The simulated probabilities of the interval $I_i^{CT}$ are more than 95\% for all situations, while the average length is larger than that of the direct confidence interval. The truncation issue of REML estimates might affect this negative performance.
Specifically, the length is about 4.6 times that of the direct confidence interval, as reported for a leverage value of 0.07 and $B=0.9$.
In contrast, $I_i^{YL}$ and $I_{i}^{NAS}$ show more than 95\% simulated coverage probabilities, and the lengths are always less than that of the direct confidence interval.
Moreover, the average length of $I_{i}^{NAS}$ is less than that of $I_{i}^{YL}$. $I_{i}^{NAS}$ is shown to improve $I_{i}^{YL}$ in terms of length.

\subsection{ Study 2: unbalanced case in which condition (\ref{cond.YL}) does not hold}
We considered the same situation as in Study 1 except for the patterns of leverage, sampling variances $D_i$, and the unknown model variance $A$. In this study, we compared the performances in an unbalanced case when condition (\ref{cond.YL}) does not hold, unlike in Study 1.
The covariate $x_{i2}$ is generated once independently from the uniform distribution $U(0,0.5)$ for the first 14 areas, and the final one is generated once from another uniform distribution, $U(0.5,1)$, and are then fixed such that the condition does not hold.
Thus, the condition is satisfied by all but the final area, with a leverage of 0.64. 
In order to simultaneously investigate the effect of the shrinkage factor $B_i$ in an unbalanced case, we designed two different unbalanced cases where the $D_i$ patterns are (a)\{0.2,0.4,0.5,0.6,2\} and (b) \{2,4,5,6,20\}, with $A$=0.1, 1 for each pattern. We assumed five groups of areas where $D_i$ is the same within each group.
Pattern (a) is the same as Pattern (b) in Datta et al. (2005), and Pattern (b) with $A$=1 shows the same $B_i$ values in Pattern (a) with $A=0.1$. Thus, $B_i$ ranges from 0.47 to 0.9. We also consider Pattern (c)$\{2,0.6,0.5,0.4,0.2\}$ with $A=0.1$ in order to investigate the effect of a moderate $B_i$ value with maximum leverage and a large $B_i$ with minimum leverage, arranging the pattern in the descending order of Pattern (a).
Incidentally, $I_{i}^{YL}$ is not included as a comparable competitor because of the failure condition required.

\begin{table}[h]
\small
\centering
\caption{Simulated coverage probabilities and average length (in parentheses) in an unbalanced case such that condition (\ref{cond.YL}) does not hold with 95\% nominal coverage}
\label{unb.Dib}
\begin{tabular}{c|c|c|cccccc}  \hline
Pattern&B & Leverage & Cox.Re & T.Re&NAS & CT.Re & Direct \\
  \hline
(a) & 0.47 & 0.07 & 78.04 & 98.04 & 96.83 & 99.2 & 94.81 \\ 
  &  &  & ( 1.05 ) & ( 1.68 ) & ( 1.57 ) & ( 5.01 ) & ( 1.75 ) \\ \cline{2-8}
  & 0.9 & 0.64 & 53.51 & 93.74 & 96.66 & 94.47 & 94.84 \\ 
  & & & ( 1.52 ) & ( 3.63 ) & ( 4.15 ) & ( 3.73 ) & ( 5.54 ) \\ \hline
  (b)& 0.47 & 0.07 & 72.21 & 97.78 & 96.47 & 99.06 & 95.28 \\ 
 &  &  & ( 3.14 ) & ( 5.33 ) & ( 4.98 ) & ( 21.2 ) & ( 5.54 ) \\ \cline{2-8}
 & 0.9 & 0.64 & 50.69 & 93.21 & 96.23 & 94.04 & 94.41 \\ 
 & & & ( 4.64 ) & ( 11.45 ) & ( 13.13 ) & ( 11.79 ) & ( 17.53 ) \\ \hline
 (c)& 0.47 & 0.64 & 72.1 & 97.74 & 95.75 & 98.33 & 95.12 \\ 
  & & & ( 1.04 ) & ( 2.0 ) & ( 1.73 ) & ( 2.3 ) & ( 1.75 ) \\ \cline{2-8}
  & 0.9 & 0.07 & 73.66 & 86.77 & 98.86 & 99.73 & 95.1 \\ 
  & &  & ( 1.51 ) & ( 1.91 ) & ( 2.83 ) & ( 4.3 ) & ( 5.54 ) \\ 
   \hline
\end{tabular}
\end{table}

The result for Study 2 is displayed in Table \ref{unb.Dib}. As in Study 1, this table reports simulated coverage probabilities with average length in parentheses for a nominal coverage of 95\% as well as minimum and maximum $B_i$ and leverage values, that is, respectively 0.47 and 0.9 for $B_i$ and 0.07 and 0.64 for leverage.
$I_i^{Cox}$ also reveals under-coverage for all situations, as in Study 1. It is noteworthy that the simulated coverage probability dramatically goes down to about 50\% in combinations with ($B_i$,Leverage)=(0.9,0.64) for both patterns (a) and (b) from this table.
$I_{i}^T$ also provides under-coverage results for large $B_i$ values. In particular, it has a considerable under-coverage problem, as much as 86.77\% despite a 95\% nominal coverage, in combination with ($B_i$,Leverage)=(0.9,0.07). The loss of Property (i) might also affect these under-coverage results.
Similar to the result in Study 1, the simulated probability of $I_{i}^{CT}$ is reported to be more than 95\% except for ($Bi$,Leverage)=(0.9,0.64), while the average lengths are mostly larger than that of the direct confidence interval.
A remarkable result is that the average interval length is about 3.8 times that of $I_i^{D}$ for ($B_i$,Leverage)=(0.47,0.07) in Pattern (b).
In contrast, $I_{i}^{NAS}$ maintains more than 95\% simulated coverage probabilities, and the length is always less than that of the direct confidence interval.
Furthermore, this study also shows that $I_{i}^{NAS}$ improves $I_{i}^{YL}$ in terms of the leverage condition.

\section{Conclusion and Discussion}

This study proposed a second-order efficient empirical Bayes confidence interval based on the non-area-specific adjustment factor under the Fay--Herriot model. Our method simultaneously provides five desired properties while the interval $I_i^{YL}$ satisfies three properties.
Additionally, the overall results from the simulation showed that our confidence interval, $I_{i}^{NAS}$, is superior to other confidence intervals, including $I_{i}^{YL}$, in terms of coverage probability and length.
Moreover, we also studied the simulated elapsed time with regard to computer burden.
Accordingly, we considered three values of $m$, $m=10,100,500$, to compare the calculation time of $I_{i}^{YL}$ with those of $I_{i}^{Cox}$, $I_{i}^{NAS}$, and $I_{i}^{D}$. In this study, we set $x_i^{\prime}\beta=\mu=0$ with $p=1$ and generated a dataset of $y_i=v_i+e_i$ with $A=1$ and $D_i=D=1$. Simulated calculation time was measured as the average of five implementations in R 3.2.0 with the Fisher scoring method and no parallel implementation for $I_{i}^{YL}$.
We found that the calculation time of $I_{i}^{YL}$ rises rapidly as $m$ increases. In particular, for $m=500$, it takes about 1000 seconds, about 1000 times as long as that of $I_{i}^{NAS}$, which takes about 1 second. We expect that much longer calculation times will be required for $I_{i}^{YL}$ with even larger $m$ values.

\section*{ Acknowledgments} The authors thank Professor Partha Lahiri at University of Maryland for reading an earlier draft of the paper and making constructive comments. This research was supported by JSPS Grant-in-Aid for Research Activity Start-up No. 26880011.

\appendix
\def\thesection{Appendix.\Alph{section}}

\section{Proof of Theorems}
\subsection{Theorem 1}

Similar to the proof of Theorem 1 in Diao et al. (2014), we have for large $m$, under regularity conditions,
\begin{align}
P[\theta_i \in I^{CT}(\hat{A}_i)]=1-\alpha+\frac{z_{\alpha/2}\phi(z_{\alpha/2})}{g_{1i}(A)}\left \{k_i+c_i-\frac{2g_{3i}(A)D_i^2}{g_{1i}(A)(A+D_i)}\right\}+O(m^{-3/2}),\label{diao}
\end{align}
where $k_i=E[g_{1i}(\hat{A}_i)]-g_{1i}(A)+g_{3i}(A)$ and $c_i=s_i^2-g_{1i}(A)-g_{2i}(A)$. Note that Diao et al. (2014) considered the existing variance estimation method for $A$.

By using the general adjusted maximum likelihood method and Theorem 1 in Yoshimori and Lahiri (2014b), we have for large $m$, under regularity conditions,
\begin{align}
E[g_{1i}(\hat{A}_i)-g_{1i}(A)]=&B_i^2\frac{2}{tr[V^{-2}]}\frac{\partial  \log \tilde{L}_{i}(A)}{\partial A}-g_{3i}(A)+o(m^{-1}).\label{g.bias}
\end{align}

Following (\ref{diao}) and (\ref{g.bias}), we have
\begin{align}
(\ref{diao})=&1-\alpha+\frac{z_{\alpha/2}\phi(z_{\alpha/2})}{g_{1i}(A)}\left \{B_i^2\frac{2}{tr[V^{-2}]}\frac{\partial  \log \tilde{L}_{i}(A)}{\partial A}+c_i-\frac{2g_{3i}(A)D_i^2}{g_{1i}(A)(A+D_i)}\right\}+O(m^{-3/2}),\notag\\
=&1-\alpha+\frac{z_{\alpha/2}\phi(z_{\alpha/2})g_{3i}(A)}{g_{1i}(A)}\left \{(A+D_i)\frac{\partial  \log \tilde{L}_{i}(A)}{\partial A}+c_i^*-\frac{2D_i}{A}\right\}+O(m^{-3/2}).\label{con}
\end{align}
where $c_i^*=c_i/g_{3i}(A)$.

With the terms on the right-hand side of (\ref{con}) vanishing, a second-order empirical Bayes confidence interval appears:
\begin{align}
(A+D_i)\frac{\partial  \log \tilde{L}_{i}(A)}{\partial A}+c_i^*-\frac{2D_i}{A}=0. \label{th.1p}
\end{align}

Thus, (\ref{ad.diff2}) follows from solving (\ref{th.1p}).

\subsection{Theorem 2 (1) and 3 (3)}

We first prove Theorem 2(1) for the existence of $\hat{A}^{NAS}$ with $A>0$. 

From the result of $\tilde{L}_{NAS}(A)= A^{\frac{(1+z^2)}{4}}$, we have
$$\tilde{L}_{NAS}(A) L_{RE}(A)\Big{|}_{A=0}=0, \ \ {\rm and} \ \ \tilde{L}_{NAS}(A) L_{RE}(A)\Big{|}_{A>0}>0.$$

Hence, it suffices to show that, for large $A$,
\begin{align*}
\tilde{L}_{NAS}(A) L_{RE}(A)=o(1).
\end{align*}

We also have, for large $A$,
\begin{align*}
\tilde{L}_{NAS}(A) L_{RE}(A)<&CA^{\frac{(1+z^2)}{4}}(A+\sup_{i\geq 1} D_{i})^{\frac{p}{2}}| X^{\prime}X |^{-\frac{1}{2}}(A+\inf_{i\geq 1}D_{i})^{-\frac{m}{2}}\\
=&O(A^{\frac{(1+z^2)}{4}-\frac{m-p}{2}}), 
\end{align*}
where $C$ denotes a generic positive constant.

Thus, the condition $m>p+\frac{1+z^2}{2}$ proves the existence of $\hat{A}^{NAS}$ for $A>0$

Similarly,
 $\tilde{L}_{i,ad}^{(c)}(A)=A^{\frac{1+z^2}{4}}(A+D_i)^{\frac{7-z^2}{4}}$ results in Theorem 3 (3).

\subsection{Theorem 2(2)}

We prove the uniqueness of $\hat{A}_{NAS}$ on $A>0$ in a balanced case.

In a balanced case, we have 
\begin{align*}
\frac{\partial \log L(A)}{\partial A}&+\frac{\partial \log \tilde{L}_{NAS}(A)}{\partial A}\\
&=\frac{1}{2(A+D)^2}\left[y^{\prime}My-(m-p)(A+D) +\frac{(1+z^2)(A+D)^2}{2A}\right].
\end{align*}
where $M=\{I_m-X(X^{\prime}X)^{-1}X^{\prime}\}$. 

Hence, our estimate $\hat{A}^{NAS}$ is obtained as a solution of
\begin{align*}
-2\left\{m-p-\frac{1+z^2}{2}\right\}A^2+2\{y^{\prime}My-(m-p-1-z^2)D\}A+(1+z^2)D^2=0.\label{uniq}
\end{align*}

Thus, Theorem 2(2) obtained from quadratic formula under $m>p+\frac{1+z^2}{2}$.

\section{Proof of $g_{1i}(A)+g_{2i}(A)<D_i$} 

From the regularity conditions R1 and R3, we have 
\begin{align*}
\frac{x_i^{\prime}(X^{\prime}V^{-1}X)^{-1}x_i}{(A+D_i)}>0.
\end{align*}

Hence, the proof follows from 
\begin{align*}
g_{1i}(A)+g_{2i}(A)=&D_i\left\{
1-\frac{D_i}{A+D_i}+\frac{D_i}{(A+D_i)}\frac{x_i^{\prime}(X^{\prime}V^{-1}X)^{-1}x_i}{(A+D_i)}
\right\}\\
=&D_i\left[1-\frac{D_i}{A+D_i}\left\{\frac{x_i^{\prime}(X^{\prime}V^{-1}X)^{-1}x_i}{(A+D_i)}+1\right\}\right]<D_i.
\end{align*}

\end{document}